\magnification=\magstep1
\hsize=6.5truein
\hoffset=0.0truein
\baselineskip 1.3\normalbaselineskip

\tolerance=10000
\def\sqr{$\vcenter{\hrule height .3mm
\hbox {\vrule width .3mm height 2mm \kern 2mm
\vrule width .3mm} \hrule height .3mm}$}
\input amssym.def
\input amssym.tex
\def \Proof{\noindent {\bf Proof.\ \ }}
\def \R{{\Bbb R}}
\def \E{{\Bbb E}}

\def \C{{\Bbb C}}

\def \seq#1#2{#1_1,\dots,#1_#2}
\def \sm#1#2{\sum_{#1=1}^#2}
\def \ge{\geqslant}
\def \le{\leqslant}

\def \dim{\mathop{\rm dim}}
\def \e{\epsilon}
\def \G{\Gamma}
\def \g{\gamma}
\def \d{\delta}
\def \a{\alpha}
\def \b{\beta}

\def \sp#1{\langle#1\rangle}
\def \ra{\rightarrow}
\def \Z{{\Bbb Z}}
\def \C{{\Bbb C}}
\def \F{{\Bbb F}}
\def \ol{\overline}
\def \hf{{\hat f}}
\def \tr{\mathop{\rm tr}}
\def \ca{{\cal A}}
\def \cb{{\cal B}}

\centerline{\bf QUASIRANDOM GROUPS}
\bigskip
\centerline {W. T. Gowers}
\bigskip

\noindent {\bf Abstract.} {\sl Babai and S\'os have asked whether
there exists a constant $c>0$ such that every finite group $G$
has a product-free subset of size at least $c|G|$: that is, a
subset $X$ that does not contain three elements $x$, $y$ and
$z$ with $xy=z$. In this paper we show that the answer is no.
Moreover, we give a simple sufficient condition for a group
not to have any large product-free subset.}
\bigskip

\noindent {\bf \S 1. Introduction.} 
\bigskip

The starting point for this paper is a well-known result of Erd\H os,
which states that for every $n$-element subset $X$ of $\Z$ there is a
subset $Y\subset X$ of size at least $n/3$ that is {\it sum-free}, in
the sense that if $y_1$ and $y_2$ belong to $Y$ then $y_1+y_2$ does
not belong to $Y$. The proof is so simple that it can be given in full
here. First, choose a prime $p$ such that $X$ lives in the interval
$[-p/3,p/3]$. A subset $Y\subset X$ is then sum-free if and only if it
is sum-free mod $p$.  But if $r$ is any integer not congruent to $0$
mod $p$, then $Y$ is sum-free mod $p$ if and only if $rY$ is sum-free
mod $p$. Moreover, a simple averaging argument shows that one can find
$r$ such that at least a third of the elements of $rX$ lie in the
interval $[p/3,2p/3]$ mod $p$. Therefore, $X$ has a subset $Y$ of size
at least $n/3$ such that $rY$, and hence $Y$, is sum-free.

Using the classification of Abelian groups it is easy to see that the
same result holds if $X$ is a subset of an Abelian group, but the
situation for non-Abelian groups is less clear.  In 1985, Babai and
S\'os [2] noted that if $H$ is a subgroup of $G$ of index $k$, then any
non-trivial coset of $H$ is product-free. From the classification of
finite simple groups it can be shown that every finite simple group of
order $n$ has a subgroup of index at most $Cn^{3/7}$ and hence a
product-free set of size at least $cn^{4/7}$. Combining that with the
fact that a product-free subset of a quotient of $G$ lifts to a
product-free subset of $G$, one can deduce the same result for all
finite groups. In 1997, Kedlaya [11] (see also [12]) improved this bound
to $cn^{11/14}$ by showing that if $H$ has index $k$ then one can find 
a union of $ck^{1/2}$ cosets of $H$, a large subset of which is product 
free.

In the other direction, nothing much was known. Indeed, Babai and S\'os
asked whether the lower bound could be improved to $cn$ for some
positive constant $c$, and Kedlaya repeated the question, while also
asking the weaker question of whether, for every $\e>0$, one can
obtain a bound of $c(\e)n^{1-\e}$.  This paper answers these questions
in the negative, by showing that, for sufficiently large $q$, the
group PSL$_2(q)$ has no product-free subset of size $Cn^{8/9}$, where
$n$ is the order of PSL$_2(q)$. In fact, we prove the stronger result
that if $A$, $B$ and $C$ are three subsets of PSL$_2(q)$ of size at
least $Cn^{8/9}$, then there is a triple $(a,b,c)\in A\times B\times
C$ such that $ab=c$.

The proof has three stages. First, we briefly review some facts about
quasirandom bipartite graphs and quasirandom subsets of groups --
detailed proofs of most of these can be found elsewhere, and we give
simple proofs of those that cannot. Secondly, we prove that the
``bipartite Cayley graph'' associated with PSL$_2(q)$ and one of the
three sets under consideration is quasirandom. Finally, we show that
this quasirandomness immediately implies our result.

Having proved this theorem, we step back and look at what we have done
from a more abstract point of view. The property of PSL$_2(q)$ that
makes it suitable for results of this kind is that it has no
non-trivial irreducible representations of low dimension. This
property has been used in a similar way before: it is due to Sarnak
and Xue [16]. It was also used in [7] to prove that the famous Ramanujan
graphs of Lubotsky, Phillips and Sarnak [14] are expanders (this is a
weaker result than that of [14] but the proof is much easier), and it
has recently been used by Bourgain and Gamburd [4] to show the same
for certain other Cayley graphs.

Our main result is rather easier than theirs. However, this very fact
may make it useful to readers who do not have a background in
representation theory and who would like to see how information about
representations can be used. If a group has no non-trivial
low-dimensional representations, it seems appropriate to call it {\it
quasirandom} since, as we show later in the paper, this property is
equivalent to several other properties, some of which state that
certain associated graphs are quasirandom. Once we have stated and
proved various equivalences of this kind, we prove some further
results. The first of these is a partial converse to our main theorem: 
if a finite group $G$ contains no large product-free subset, then it is
quasirandom. The reason this is a ``partial'' converse is that the
bounds we obtain are not very good: for most of the results in the
paper there is a power-type dependence of one constant on another,
but for this one it is exponential/logarithmic.

Section 4 ends with another weak equivalence. It is easy to prove
that a group is not quasirandom if it has a non-trivial quotient
that is either Abelian or of small order. We show that, in the
absence of these obvious obstructions, a group $G$ is quasirandom.
In particular, non-Abelian finite simple groups are quasirandom.
Again, we obtain exponential/logarithmic bounds, but for this result
it is unavoidable because the dimension of the smallest non-trivial 
representation is a power of $n$ for some finite simple groups
and logarithmic in $n$ for others.

In Section 5 we prove a generalization of the main theorem to
more complicated sets of equations. The theorem itself allows one to
place $a$, $b$ and $ab$ into specified dense subsets of a quasirandom
group. It turns out that one can do the same with more variables: for
example, the next case says that $a$, $b$, $c$, $ab$, $bc$, $ac$ and 
$abc$ can be placed into specified sets. 

The final section of this paper collects together some open problems
that have arisen during the paper, and adds a few more. 
\bigskip

\noindent {\bf \S 2. Quasirandom graphs and sets.}
\bigskip

As promised, let us briefly review some of the standard theory of
quasirandomness, concentrating in particular on the definitions of a
quasirandom graph, a quasirandom bipartite graph and of a quasirandom
subset of an Abelian group. The first few results of this section will
not be used later, so we shall not give their proofs. However, they
put the later results into their proper context.

The notion of a quasirandom graph was introduced by Chung, Graham and
Wilson [6], though a similar notion (of so-called ``jumbled''
graphs) had been defined by Thomason [17]. If $x$ is a vertex in a
graph, we shall write $N_x$ for its neighbourhood. The {\it adjacency
matrix} $A$ of a graph $G$ is defined by $A(x,y)=1$ if $xy$ is an edge
of $G$ and $A(x,y)=0$ otherwise.
\bigskip

\noindent {\bf Theorem {2.1}.} {\sl Let $G$ be a graph with $n$ vertices 
and density $p$. Then the following statements are polynomially equivalent, 
in the sense that if one statement holds for a constant $c$, then 
all others hold with constants that are bounded above by a positive
power of $c$.

(i) $\sum_{x,y\in V(G)}|N_x\cap N_y|^2\le (p^4+c_1)n^4$.

(ii) The number of labelled 4-cycles in $G$ is at most $(p^4+c_1)n^4$.

(iii) For any two subsets $A,B\subset V(G)$ the number of 
pairs $(x,y)\in A\times B$ such that $xy\in E(G)$ differs
from $p|A||B|$ by at most $c_2n^2$.

(iv) The second largest modulus of an eigenvalue of the adjacency
matrix of $G$ is at most $c_3n$.}
\bigskip

A graph that satisfies one, and hence all, of these properties for a
small $c$ is called {\it quasirandom}. If one wishes to be more
precise, then one can say that $G$ is $c$-quasirandom if it satisfies
property (i) (or equivalently (ii)) with constant $c_1=c$. A random
graph with edge probability $p$ is almost always quasirandom with
small $c$, and quasirandom graphs have many properties that random
graphs have. In particular, if $H$ is any fixed small graph, and
$\phi$ is a random map from $V(H)$ to $V(G)$, then the probability
that $\phi(x)\phi(y)$ is an edge of $G$ whenever $xy$ is an edge of
$H$ (in which case $\phi$ is a {\it homomorphism}) is roughly
what one would expect, namely $p^{|E(H)|}$, and the probability that
in addition no non-edge of $H$ maps to an edge of $G$ (in which case
$\phi$ is an {\it isomorphic embedding}) is roughly
$p^{|E(H)|}(1-p)^{{|V(H)|\choose 2}-|E(H)|}$. 

A quasirandom bipartite graph is like a quasirandom graph but with
some obvious modifications. As above, we state a theorem that serves
as a definition as well. 
\bigskip

\noindent {\bf Theorem {2.2}.} {\sl Let $G$ be a bipartite graph 
with vertex sets $X$ and $Y$ and $p|X||Y|$ edges. Then the following
statements are polynomially equivalent.

(i) $\sum_{x,x'\in X}|N_x\cap N_{x'}|^2\le (p^4+c_1)|X|^2|Y|^2$.

(i) $\sum_{y,y'\in Y}|N_y\cap N_{y'}|^2\le (p^4+c_1)|X|^2|Y|^2$

(ii) The number of labelled 4-cycles that start in $X$ is at most 
$(p^4+c_1)|X|^2|Y|^2$.

(iv) For any two subsets $A\subset X$ and $B\subset Y$ the number of 
pairs $(x,y)\in A\times B$ such that $xy\in E(G)$ differs
from $p|A||B|$ by at most $c_2|X||Y|$.}
\bigskip

\noindent We call a bipartite graph $c$-{\it quasirandom} if it 
satisfies condition (i) (and therefore the exactly equivalent
conditions (ii) and (iii)) with constant $c_1=c$. 

Note that we have not given an eigenvalue condition. This
is because the bipartite adjacency matrix (that is, the obvious
01-function defined on $X\times Y$ as opposed to $(X\cup Y)^2$)
is not symmetric. However, as we shall see later, there is a natural
analogue of this condition.

To continue our quick survey of known results, let us define
quasirandom subsets of Abelian groups. This is a straightforward
generalization of a definition of Chung and Graham [5] for the case
$\Z/p\Z$. Again, we present it as a theorem rather than a
definition. Recall that if $G$ is an Abelian group, $f$ is a function
from $G$ to $\C$ and $\g:G\ra\C$ is a character of $G$, then the
Fourier transform of $f$, evaluated at $\g$, is the number
$\hf(\g)=|G|^{-1}\sum_{g\in G}f(g)\overline{\g(g)}$. If $f_1$ and
$f_2$ are two functions defined on $G$, then their {\it convolution}
$f_1*f_2$ is defined by $f_1*f_2(g)=\sum_{x+y=g}f_1(x)f_2(y)$. If $A$
is a subset of $G$ we shall use the letter $A$ also for the
characteristic function of $A$. That is, $A(x)=1$ if $x\in A$ and 0
otherwise.
\bigskip

\noindent {\bf Theorem {2.3}.} {\sl Let $G$ be an Abelian group of
order $n$ and let $A\subset G$ be a set of size $pn$. Then the following 
are equivalent.

(i) $\sum_{g\in G}|A\cap(A+g)|^2\le (p^4+c_1)n^3$.

(ii) There are at most $(p^4+c_1)n^3$ solutions in $A$ of the equation
$x+y=z+w$.

(iii) $\sum_{g\in G}|A*A(g)|^2\le (p^4+c_1)n^3$.

(iv) For every subset $B\subset G$, $\sum_{g\in G}|A*B(g)|^2\le 
n^{-1}|A|^2|B|^2+c_2n^3$.

(v) The graph with vertex set $G$ and with $x$ joined to $y$ if
and only if $x+y\in A$ is $c_1$-quasirandom.

(vi) The bipartite graph with two copies of $G$ as its vertex sets
and with $x$ joined to $y$ if and only if $y-x\in A$ is $c_1$-quasirandom.

(vii) $|{\hat A}(\g)|\le c_3n$ for all non-trivial characters $\g$.}
\bigskip

It is often convenient to replace Theorems 2.2 and 2.3 with
``functional'' or ``analytic'' versions, as follows.

\noindent {\bf Theorem {2.4}.} {\sl Let $X$ and $Y$ be two finite
sets and let $f:X\times Y\ra\C$ be a function that takes values
of modulus at most 1. Then the following properties of $f$ are
polynomially equivalent.

(i) $\sum_{x,x'\in X}\sum_{y,y'\in Y}f(x,y)\overline{f(x,y')}
\overline{f(x',y)}f(x',y')\le c_1|X|^2|Y|^2$.

(ii) For any two functions $u:X\ra\C$ and $v:Y\ra\C$ taking values
of modulus at most~1, 
$$\Bigl|\sum_{x,y}f(x,y)u(x)v(y)\Bigr|\le c_2|X||Y|.$$
 
(iii) For any two sets $A\subset X$ and $B\subset Y$,
$$\Bigl|\sum_{x\in A}\sum_{y\in B}f(x,y)\Bigr|\le c_3|X||Y|.$$}

\noindent A function $f$ with one, and hence all three, of the above
properties is called {\it quasirandom}. More precisely, we call it
$c$-quasirandom if property (i) holds with constant $c$.

Theorem 2.4 is closely related to Theorem 2.2. Indeed,
if $G$ is a bipartite graph with vertex sets $X$ and $Y$ and 
density $p$, then $G$ is quasirandom if and only if the function
$f(x,y)=G(x,y)-p$ is quasirandom, where we have written $G$ for
the characteristic function of the graph as well (so $f(x,y)$ is 
$1-p$ if $(x,y)$ is an edge and $-p$ otherwise). This is particularly
easy to show if $G$ is {\it regular}, in the sense that every vertex 
in $X$ has degree $p|Y|$ and every vertex in $Y$ has degree $p|X|$.
Then a quick calculation shows that $G$ is $c$-quasirandom if and 
only if $f$ is $c$-quasirandom.

Now let us give a functional version of Theorem 2.3. Instead of
trying to give as many equivalences as possible, we shall restrict
our attention to ones that will be of interest later (in Section 4,
when we come to define quasirandom groups). These apply to subsets
of an arbitrary group. They are not deep equivalences, as one might
suspect from the fact that they all hold with the same constant.
\bigskip

\noindent {\bf Theorem {2.5}.} {\sl Let $G$ be a group of
order $n$ and let $f:G\ra\C$ be a function taking values of modulus
at most 1. Then the following are exactly equivalent. 

(i) $\sum_{x\in G}\Bigl|\sum_{y\in G}f(x)\overline{f(yx)}\Bigr|^2
\le cn^3$.

(ii) $\sum_{ab^{-1}=cd^{-1}}f(a)\overline{f(b)}\overline{f(c)}f(d)\le cn^3$.

(iii) The function $F(x,y)=f(xy^{-1})$ is a $c$-quasirandom function
on $G\times G$.}

\Proof  To see that (i) and (ii) are equivalent, note that the sum
on the left-hand side of (i) is equal to
$$\sum_{x,y,z\in G}f(x)\overline{f(yx)}\overline{f(z)}f(yz).$$
The result now follows from the obvious one-to-one correspondence between
quadruples $(a,b,c,d)$ such that $ab^{-1}=cd^{-1}$ and quadruples of the form
$(x,yx,z,yz)$.

To see that (ii) and (iii) are equivalent, note that
$$\sum_{x,x'}\sum_{y,y'}F(x,y)\overline{F(x,y')}\overline{F(x',y)}F(x',y')
=\sum_{x,x'}\sum_{y,y'}f(xy^{-1})\overline{f(xy'^{-1})}
\overline{f(x'y^{-1})}f(x'y'^{-1})\ .$$
Now for each $x,x',y$ and $y'$ we have 
$(xy^{-1})(x'y^{-1})^{-1}=(xy'^{-1})(x'y'^{-1})^{-1}$. In the
other direction, if $ab^{-1}=cd^{-1}$ and $g$ is any group element,
then let $y=g$, $x=ag$, $y'=c^{-1}ag$ and $x'=dc^{-1}ag=bg$. Then
$xy^{-1}=a$, $x'y^{-1}=b$, $xy'^{-1}=c$ and $x'y'^{-1}=d$. This gives
us an $n$-to-one correspondence between quadruples 
$(xy^{-1},x'y^{-1},xy'^{-1},x'y'^{-1})$ and quadruples $(a,b,c,d)$
such that $ab^{-1}=cd^{-1}$, which proves that (ii) holds if and
only if
$$\sum_{x,x'}\sum_{y,y'}F(x,y)\overline{F(x,y')}\overline{F(x',y)}F(x',y')
\le cn^4,$$
that is, if and only if (iii) holds. \hfill $\square$ \bigskip

If these properties hold (as well as the hypotheses of the theorem)
then we shall say that $f$ is $c$-{\it quasirandom}. For more details
about quasirandom graphs, sets and functions, including proofs of most
of the previous results, the reader is referred to the early sections
of [9]. (This is by no means the only reference, but is chosen because
the presentation there harmonizes well with the presentation in this
paper.)
\bigskip

Let us now return to the question of a ``spectral theory'' for
bipartite graphs. For an ordinary graph $G$, one observes that the
adjacency matrix is symmetric and can therefore be decomposed as $\sm
i n \lambda_i u_i\otimes u_i$ for some orthonormal basis $(u_i)$ of
eigenvectors, with $\lambda_i$ the eigenvalue corresponding to
$u_i$. (Here we write $u\otimes v$ for the matrix that takes the value
$u(x)v(y)$ at $(x,y)$. If $v$ and $w$ are elements of inner product
spaces $V$ and $W$, then we write $w\otimes v$ for the linear map from
$V$ to $W$ defined by $x\mapsto\sp{x,v}w$. Notice that these two
definitions are consistent.) For a bipartite graph, the adjacency
matrix is no longer symmetric, so this result is no longer
true. However, what we can do instead is decompose it as a sum $\sm i
n\lambda_i u_i\otimes v_i$, where $(u_i)$ and $(v_i)$ are {\it two}
orthonormal bases. This is called the {\it singular value
decomposition} of the matrix, which was discovered in the late 19th
century and is important in numerical analysis. For the convenience of
the reader, we give a proof that it always exists (in the real case).

\proclaim Theorem {2.6}. Let $\a$ be any linear map from a real
inner product space $V$ to a real inner product space $W$. 
Then $\a$ has a decomposition of the form $\sm i k\lambda_iw_i\otimes v_i$,
where the sequences $(w_i)$ and $(v_i)$ are orthonormal in $W$
and $V$, respectively, each $\lambda_i$ is non-negative, and $k$ is 
the smaller of $\dim V$ and $\dim W$.

\Proof  To begin, let $v$ be a non-zero vector such that 
$\|\a v\|/\|v\|$ is maximized. (For this proof, $\|.\|$ is the
standard Euclidean norm and $\sp{,}$ the standard inner product,
either on $\R^m$ or $\R^n$.) Now suppose that $w$ is any vector
orthogonal to $v$ and let $\d$ be a small real number. Then $\|\a(v+\d
w)\|^2=\|\a v\|^2+2\d\sp{\a v,\a w}+o(\d)$, and $\|v+\d
w\|^2=\|v\|^2+o(\d)$. It follows that $\sp{\a v,\a w}=0$, since
otherwise we could pick a small $\d$ with the same sign as $\sp{\a
v,\a w}$ and we would find that $\|\a(v+\d w)\|/\|v+\d w\|$ was bigger
than $\|\a v\|/\|v\|$.

Let $X$ and $Y$ be the subspaces of $\R^n$ and $\R^m$ orthogonal to
$v$ and $\a v$, respectively. They can be given orthonormal bases, and
$\a$ maps everything in $X$ to $Y$. Let $\b$ be the restriction of
$\a$ to $X$. By induction, $\b$ has a decomposition of the required
form. That is, we can write $\b=\sum_{i=2}^k\lambda_iw_i\otimes v_i$
with $v_i\in X$ and $w_i\in Y$. Now set $v_1=v/\|v\|$, $w_1=\a v/\|\a v\|
=\a v_1/\|\a v_1\|$ and $\lambda_1=\|\a v_1\|$. Then $\a v_1=\lambda_1w_1$,
from which it follows that $\a=\sm i k\lambda_iw_i\otimes v_i$, as
required. \hfill $\square$ \bigskip

This theorem is of course equivalent to a very similar statement 
about matrices, and indeed that is how we shall apply it.

The fact that singular values are the correct analogue of eigenvalues
for bipartite graphs has been realized before. See for example [3]. 
The next two results illustrate the connection very clearly.

\proclaim Lemma {2.7}. Let $G$ be a bipartite graph with vertex 
sets $X$ and $Y$ and identify $G$ with its bipartite adjacency matrix 
$\sm i k\lambda_iw_i\otimes v_i$, where $(v_i)$ and $(w_i)$ are 
orthonormal sequences. Then $\sum_i\lambda_i^2$ is the number of 
edges in $G$ and $\sum_i\lambda_i^4$ is the number of labelled
4-cycles that start in $X$. 

\Proof  The number of edges in $G$ is $\tr(G^TG)$. But
$G^T$ is $\sum_i\lambda_iv_i\otimes w_i$. It is easy to verify
that $(v_i\otimes w_i)(w_j\otimes v_j)=v_i\otimes v_j$. But
$\tr(v_i\otimes v_j)=1$ if $i=j$ and $0$ otherwise, so the
first statement of the lemma follows. 

The second part is similar. The number of labelled 4-cycles that
start in $X$ is $\tr(G^TGG^TG)$. If we expand $G$ and $G^T$ then
once again the only terms that survive are those that use a single
$i$. But in this case we have four terms, so the answer is
$\sum_i\lambda_i^4$.~\hfill~$\square$ \bigskip
 
The next result gives a further condition that is 
equivalent to quasirandomness for regular bipartite graphs.
\bigskip

\noindent {\bf Theorem {2.8}.} {\sl Let $G$ be a regular bipartite graph 
with vertex sets $X$ and $Y$, $p|X||Y|$ edges and identify $G$ with 
its bipartite adjacency matrix. Then the following are polynomially equivalent.

(i) $G$ is $c_1$-quasirandom.

(ii) The maximum of $\|Gf\|/\|f\|$ over all non-zero functions 
$f$ such that $\sum_{x\in X}f(x)=0$ is at most $c_2|X|^{1/2}|Y|^{1/2}$.}

\Proof  By Theorem 2.6 we can write $G=\sm i k\lambda_iw_i\otimes v_i$ 
for orthonormal sequences $(v_i)$ and $(w_i)$. By Lemma 2.7, the number 
of labelled 4-cycles in $G$ that start in $X$ is $\sm i k\lambda_i^4$.
Suppose that the decomposition is chosen so that $u_1$ and $v_1$ are
constant functions, which implies that $\lambda_1=p|X|^{1/2}|Y|^{1/2}$.
Then, if (ii) holds, we find that 
$$\sm i k\lambda_i^4\le p^4|X|^2|Y|^2+c_2^2|X||Y|\sum_{i=2}^k\lambda_i^2\ .$$
By Lemma 2.7, $\sum_{i=2}^k\lambda_i^2\le p|X||Y|$, so this is at most
$(p^4+pc_2^2)|X|^2|Y|^2$, which establishes (i) with $c_1=pc_2^2$.

Conversely, if (i) holds, then $\sm i k\lambda_i^4\le(p^4+c_1)|X|^2|Y|^2$.
Since $\lambda_1=p|X|^{1/2}|Y|^{1/2}$, it follows that every other 
$\lambda_i$ is at most $c_1^{1/4}|X|^{1/2}|Y|^{1/2}$. Since the maximum
of these other $\lambda_i$ is precisely the maximum in (ii), we have
established (ii) with $c_2=c_1^{1/4}$. \hfill $\square$ \bigskip

The next lemma is a simple fact, but for our purposes it will be very
important. In the statement, if $G$ is a bipartite graph with vertex
sets $X$ and $Y$ of not necessarily the same size, we call it {\it
regular} if every vertex in $X$ has the same degree and every vertex
in $Y$ has the same degree.

\proclaim Lemma {2.9}. Let $G$ be a regular bipartite graph with vertex 
sets $X$ and $Y$. Let $\a$ be the linear map from $\C^X$ to $\C^Y$
derived from the bipartite adjacency matrix of $G$. (That is, if
$f:X\ra\C$ then $\a f(y)=\sum_{x\in X,xy\in E(G)}f(x)$.)  Then the set
of all functions $f:X\ra\C$ such that $\sum_{x\in X}f(x)=0$ and $\|\a
f\|/\|f\|$ is maximized forms a linear subspace of $\C^X$.

\Proof  Let us first check, using the regularity of $G$, that the
maximum of $\|\a f\|/\|f\|$ over {\it all} functions is attained
when $f$ is a constant function. Let every vertex in $X$ have 
degree $p|Y|$, so that every vertex in $Y$ has degree $p|X|$. Then,
settting $G(x,y)$ to be 1 if $xy\in E(G)$ and $0$ otherwise,
$$\eqalign{\|\a f\|^2&=\sum_y\Bigl|\sum_xf(x)G(x,y)\Bigr|^2\cr
&=\sum_{x,x'}f(x)\ol{f(x')}\sum_yG(x,y)G(x',y)\cr
&\le{1\over 2}\sum_{x,x'}\bigl(|f(x)|^2+|f(x')|^2\bigr)\sum_y
G(x,y)G(x',y)\cr
&=\sum_x|f(x)|^2\sum_{x'}\sum_yG(x,y)G(x',y)\cr
&=\sum_x|f(x)|^2p^2|X||Y|=p^2|X||Y|\|f\|^2\ .\cr}$$
It follows that $\|\a f\|/\|f\|$ never exceeds $p|X|^{1/2}|Y|^{1/2}$.
This bound is attained when $f$ is the constant function $1$:
then $\|f\|=|X|^{1/2}$, and $\|\a f\|=p|X||Y|^{1/2}$ since $\a f$ 
takes the value $p|X|$ everywhere on $Y$.

The proof of Theorem 2.6 now
tells us that the restriction of the linear map $\a$ to
the space of functions that sum to zero can be decomposed as
$\sum_{i=2}^n\lambda_iw_i\otimes v_i$.  Without loss of generality,
$\lambda_2\ge\dots\ge\lambda_n\ge 0$. Choose $k$ such that
$\lambda_2=\dots=\lambda_k>\lambda_{k+1}$ and let $X$ be the subspace
of $G^\C$ generated by $v_2,\dots,v_k$. Then the restriction of $\a$
to $X$ is $\lambda_2\sum_{i=2}^kw_i\otimes v_i$. This map is orthogonal
on to its image, so $\|\a f\|=\lambda_2\|f\|$ for every $f\in X$. 
Since $\a\Bigl(\sum_{i=2}^n\mu_iv_i\Bigr)=\sum_{i=2}^n\lambda_i\mu_iw_i$,
it is clear that $\|\a f\|<\lambda_2\|f\|$ whenever $\sum_{x\in G}f(x)=0$
and $f\notin X$. \hfill $\square$ \bigskip

\noindent {\bf \S 3. A group with no large product-free subset.}
\bigskip

In this section we give a quick proof that the density of the largest
product-free subset of the group PSL$_2(q)$ tends to zero as $q$ tends
to infinity. Recall that PSL$_2(q)$ is the 2-dimensional projective
special linear group over $\F_q$, that is, the group of all $2\times
2$ matrices over $\F_q$ with determinant 1, quotiented by the subgroup
consisting of $I$ and $-I$. It is natural to look at this family of
groups, since it is one of the simplest infinite families of finite
simple groups; simple groups themselves are natural to look at because
if $G'$ is a quotient of a group $G$, then any product-free subset of
$G'$ lifts to a product-free subset of $G$. As we have already
mentioned, our proof will depend on one basic fact about
representations of PSL$_2(q)$, which we state without proof.

\proclaim Theorem {3.1}. Every non-trivial representation of 
PSL$_2(q)$ has dimension at least $(q-1)/2$. \hfill $\square$

\noindent The proof of Theorem 3.1, due to Frobenius, is not especially 
hard, though it isn't trivial either. A nice presentation of it can be
found in [7]. To put this result in perspective, the order of
PSL$_2(q)$ is $q(q^2-1)/2$, so the lowest dimension of a non-trivial
representation is proportional to the cube root of the order of the
group. This tells us that, in a certain sense, PSL$_2(q)$ is very far
from being Abelian.

As mentioned in the introduction, we shall in fact prove a result that
is more general in several ways. First of all, we shall prove it for
{\it any} group $\Gamma$ that has no low-dimensional non-trivial
representation.  Secondly, we shall prove an ``off-diagonal'' result:
given any {\it three} large subsets $A$, $B$ and $C$ of $\Gamma$,
there is a triple $(a,b,c)\in A\times B\times C$ such that $ab=c$.  In
order to prove this, it will be convenient (though not essential) to
express the number of such triples in terms of the following bipartite
Cayley graph $G$. The two vertex sets of $G$ are copies of $\Gamma$
and $xy$ is an edge if and only if there exists $a\in A$ such that
$ax=y$. (Note that if $xy$ is an edge, it does not follow that $yx$ is
an edge -- this is why we have to consider bipartite graphs.) Then the
number of triples we are trying to count is the number of edges from
the copy of $B$ on one side of this bipartite graph to the copy of $C$
on the other. If $|\Gamma|=n$ and $r=|A|/n$, then we know from Theorem
2.2 that the number of edges between these copies of $B$ and $C$ will
be approximately $r|B||C|$ if $G$ is sufficiently quasirandom.

We shall make this argument precise later in the section. But first,
let us prove that the graph $G$ actually is quasirandom. 

\proclaim Lemma {3.2}. Let $\Gamma$ be a finite group and suppose
that $\Gamma$ has no non-trivial representation of dimension less
than $k$. Let $A$ be any subset of $\Gamma$ and let 
$G$ be the bipartite Cayley graph defined above. Let $\a$ be the 
corresponding linear map defined in the statement of Lemma 2.5. 
Let $f:\Gamma\ra\C$ be any function such that 
$\sum_{x\in\Gamma}f(x)=0$. Then $\|\a f\|/\|f\|\le(|A|n/k)^{1/2}$.

\Proof  Note first that, for any $x$ and $y$ in $\G$, there exists
$a\in A$ such that $ax=y$ if and only if $yx^{-1}\in A$. Thus, this
is another way of stating which pairs $xy$ are edges of $G$. Writing $A$
for the characteristic function of the set $A$, we now have
$$\a f(y)=\sum_xG(x,y)f(x)=\sum_xA(yx^{-1})f(x)=\sum_{uv=y}A(u)f(v)
=A*f(y)\ ,$$
where the last equality is true by the definition of the convolution
of two functions defined on an arbitrary group. That is, $\a f=A*f$.

Let $\lambda$ be the maximum of $\|\a f\|/\|f\|$ over all
functions $f$ that sum to zero, and let $X$ be the set of all functions 
$f$ that achieve this maximum. Then $X$ is a linear subspace of
$\C^\G$, by Lemma 2.9 (of course, we count 0 as belonging to
$X$). Now if we choose any $f\in X$ and any group element
$g\in\Gamma$, then the function $T_gf$, defined by $T_gf(x)=f(xg)$,
also belongs to $X$, since 
$$\a T_gf(u)=\sum_{xy=u}A(x)T_gf(y)=\sum_{xy=u}A(x)f(yg)
=\sum_{xy=ug}A(x)f(y)=\a f(ug)\ ,$$
from which it follows that $\|\a T_gf\|=\|\a f\|$. Obviously,
$\|T_gf\|=\|f\|$ as well.

Since any non-zero $f$ in $X$ is non-constant, there exists $g\in\Gamma$ such 
that $T_gf\ne f$, from which it follows that the right-regular
representation of $\Gamma$ acts non-trivially on $X$. Therefore,
the dimension of $X$ is at least $k$, by hypothesis.

It follows from Theorem 2.6 and Lemma 2.7 that $k\lambda^2$
is at most the number of edges in $G$, which is $|A|n$. That is,
$\lambda\le(|A|n/k)^{1/2}$, as stated. \hfill $\square$ \bigskip

We have shown that $G$ satisfies condition (ii) of Theorem 2.8, with 
$c_2=(|A|/kn)^{1/2}$, as stated. This may make it look as though
$G$ becomes more quasirandom as the cardinality of $A$ decreases, but
that is just an accident arising from the way the condition is
formulated. The point is that when $A$ is smaller, the graph is less 
dense, which makes it hard for $c_2$ to be small enough for condition
(iv) of Theorem 2.2 to say anything non-trivial.

Nevertheless, we have more or less proved the main result of this
paper. All that remains is to put together the results we have
stated or proved already. 

\proclaim Theorem {3.3}. Let $\G$ be a finite group with no non-trivial
representation of dimension less than $k$, let $n=|\G|$ and let
$A$, $B$ and $C$ be three subsets of $\Gamma$ such that 
$|A||B||C|>n^3/k$. Then there exist $a\in A$, $b\in B$
and $c\in C$ with $ab=c$. In particular, this is true if all of
$A$, $B$ and $C$ have size greater than $n/k^{1/3}$. Furthermore,
if $\eta>0$ and $|A||B||C|\ge n^3/\eta^2k$, then the number
of triples $(a,b,c)\in A\times B\times C$ such that $ab=c$ is
at least $(1-\eta)|A||B||C|/n$.

\Proof  Let $|A|=rn$, $|B|=sn$ and $|C|=tn$. As in the previous
lemma, let $\a$ be the linear map $f\mapsto A*f$. Let $B$ stand
for the characteristic function of the set $B$, and for each
$x\in\G$ let $f(x)=B(x)-s$. Then $\sum_xf(x)=0$, and 
$\|f\|^2=(1-s)^2|B|+s^2(n-|B|)=s(1-s)n\le sn$.

It follows from Lemma 3.2 that $\|\a f\|^2\le rn^2sn/k$.
But $A*B(y)=A*(f+s)(y)=\a f(y)+rsn$, so whenever $A*B(y)=0$ 
we have $|\a f(y)|=rsn$. It follows that the number $m$ of $y$
for which $A*B(y)=0$ satisfies the inequality 
$m(rsn)^2\le rsn^3/k$, or $m\le n/rsk$. But if
$rst>1/k$ then this is less than $tn$, which implies
that there exists $c\in C$ such that $A*B(c)\ne 0$.
Equivalently, there exist $a\in A$ and $b\in B$ such 
that $ab=c$, as claimed. 

As for the final claim, the number of triples in question is
$\sp{A*B,C}=\sp{\a f,C}+rsn|C|$. But 
$|\sp{\a f,C}|^2\le rn^2sn|C|/k=|A||B||C|n/k$,
by the Cauchy-Schwarz inequality and the estimate for 
$\|\a f\|$ obtained earlier, while $rsn|C|=|A||B||C|/n$.
The result is therefore true provided
$$|A||B||C|n/k\le \eta^2|A|^2|B|^2|C|^2/n^2\ ,$$
and this inequality follows from our assumption.
\hfill $\square$ \bigskip

Recently, Kedlaya [13] proved a sort of converse to Theorem 3.3: under
the additional hypothesis that $G$ admits a transitive action on a
reasonably large finite set, there exist sets $A$, $B$ and $C$ such
that $|A||B||C|\ge c|\G|^3/k$ and such that there do not exist
$a\in A$, $b\in B$ and $c\in C$ with $ab=c$.

Theorems 3.1 and 3.3 immediately give the following corollary, which 
is the result promised at the beginning of the section.

\proclaim Corollary {3.4}. Let $\Gamma$ be the group PSL$_2(q)$ and
let $n=|\Gamma|$. Then $\Gamma$ has no product-free subset of cardinality 
greater than $2n^{8/9}$.

\Proof  This follows from the Theorems 3.1 and 3.3, since $n=q(q^2-1)/2$ and
$k$ can be taken to be $(q-1)/2$, which is greater than $n^{1/3}/8$.
\hfill $\square$ \bigskip 

\noindent {\bf \S 4. Quasirandom groups.}
\bigskip

The property we have just used for showing that a group $\Gamma$
does not contain a large product-free set was that $\Gamma$ has
no non-trivial low-dimensional representations. From this we deduced
that every large subset of $\Gamma$ gives rise to a directed Cayley 
graph that is quasirandom. Now we shall show that these two
properties, as well as several others, are in fact equivalent.  We
shall use the word ``quasirandom'' for any group that has one, and
hence all, of these properties, but there is a limit to how seriously
this word should be taken. In particular, we do not have a model of
random groups for which we can show that almost every group is
quasirandom. (Gromov has, famously, defined a notion of random group,
by taking a set of $n$ generators and a certain number of random
relations of prescribed length. However, his groups are infinite: to
define a random finite group one would need enough relations to make
it finite, but not enough to make it trivial, or very small. This
could be a delicate matter.)

A second difference between this notion of quasirandomness and
the usual ones for graphs and subsets of groups is that we do 
not have a ``local'' characterization, where we count small
configurations of a certain kind. (For graphs and subsets of
groups these configurations are 4-cycles and quadruples
$ab^{-1}=cd^{-1}$, respectively.) Indeed, it seems quite likely
that no such characterization exists, and to see why, consider
the case of the group $S_n$. This is not quasirandom, since
$A_n$ is a subgroup of index 2, but if you choose a small number
of permutations $\seq \pi k$ at random (here $k$ should be 
thought of as an absolute constant), then they will not have
any small relations, so one will not have any ``local'' evidence
that they are not all even permutations. That is, $S_n$ appears
to be ``locally indistinguishable'' from $A_n$, which is quasirandom.

This may not be the end of the story, however, because there is
a sense in which the non-quasirandomness of $S_n$ is at least 
``polynomially detectable.'' Suppose that you are given the 
multiplication table of $S_n$, but you are given it abstractly
and not told the order in which the permutations appear. Now
suppose that you want an algorithm that will partition the
elements into even and odd permutations in polynomial time
(in $n!$). You can do it with a randomized algorithm as follows. 
Choose $k$ elements at random from the group. Then the probability 
that they all happen to be even permutations is $2^{-k}$, and it is 
known that if they are all even then they almost surely generate 
$A_n$, while if they are not all even then they almost surely
generate $S_n$. The time it takes to find the subgroup they 
generate is easily seen to be polynomial, so after a few attempts
one will almost certainly generate $A_n$ (and we will know that
we have done so, since $A_n$ is the only subgroup of $S_n$ of
index 2). For a more general discussion of algorithms to find
the irreducible representations of a group $G$, see [1].

Now let us begin the process of
proving the main result of the section, the statement that various
properties of groups are equivalent. Before we get to the statement
itself, we shall need some mostly standard lemmas.

\proclaim Lemma {4.1}. Let $S$ be the unit sphere in $\C^n$ in the
standard Euclidean norm, and let $\mu$ be the standard rotation-invariant
probability measure on $S$. Then $\int\int|\sp{v,w}|^2d\mu(v)d\mu(w)=n^{-1}$.

\Proof  The integral in question is the mean square of the inner product
of two random unit vectors. This average is clearly unaffected if we
fix one of the vectors. But if $(e_i)_{i=1}^n$ is an orthonormal basis
of $\C^n$, then $\int_S\sm i n|\sp{v,e_i}|^2d\mu(v)=\int_S1d\mu(v)=1$,
so by symmetry $\int_S|\sp{v,e_1}|^2d\mu(v)=n^{-1}$. This proves
the lemma. \hfill $\square$ \bigskip

\proclaim Lemma {4.2}. Let $\a$ be a linear map from $\C^n$ to
$\C^n$. Then $\tr(\a)=n\int_S\sp{\a v,v}d\mu$.

\Proof  Let $(e_i)_{i=1}^n$ be an orthonormal basis. Then the
trace of the matrix of $\a$ with respect to this basis, and
hence of $\a$ itself, is $\sm i n\sp{\a e_i,e_i}$. Since this
is true for any orthonormal basis, we may average over all of
them. The result follows immediately. \hfill $\square$ \bigskip

\proclaim Lemma {4.3}. Let $v_1$ and $v_2$ be two vectors in 
$\C^n$. Then $\sp{v_1,v_2}=n\int_S\sp{v_1,w}\sp{w,v_2}d\mu(w)$.

\Proof  The proof is basically the same as that of Lemma 4.2,
since for any orthonormal basis 
$\sp{v_1,v_2}=\sm i n\sp{v_1,e_i}\sp{e_i,v_2}$, and once again
we can average over all of them. \hfill $\square$ \bigskip

\proclaim Lemma {4.4}. Let $\seq v n$ be unit vectors in $\C^m$.
Then $\sum_{i,j}|\sp{v_i,v_j}|^2\ge m^{-1}n^2$.

\Proof  The trick here is to notice that 
$|\sp{v_i,v_j}|^2=\sp{v_i\otimes\ol{v_i},v_j\otimes\ol{v_j}}$,
where $v_i\otimes\ol{v_i}$ is the $m\times m$ matrix with
entries $v_i(p)\ol{v_i(q)}$, and the inner product is the
standard inner product on $\C^{m^2}$. It follows that
$$\sum_{i,j}|\sp{v_i,v_j}|^2=\Bigl\|\sm i n v_i\otimes\ol{v_i}\Bigr\|^2\ .$$
Now $\tr(v_i\otimes\ol{v_i})=1$ for each $i$, so the trace of
$\sm i n v_i\otimes\ol{v_i}$ is $n$, from which it follows that
the right hand side is at least $m^{-1}n^2$, which proves the lemma. 
\hfill $\square$ \bigskip

Note that Lemma 4.4 is sharp. Basically any sufficiently symmetric
example shows this, but one simple one is when $m|n$ and the vectors
$v_i$ consist of $n/m$ copies of some orthonormal basis. Lemma 4.1 
proves that the result is sharp for a ``continuous set'' of vectors.
Given a set for which the lemma is sharp, the proof above shows that
$\sm i n v_i\otimes\ol{v_i}$ is $n/m$ times the identity matrix.
That is, the vectors $v_i$ give us a representation of the identity,
which is a well-known way of saying that they are nicely distributed
round the unit sphere.

With these lemmas in place, we are ready for our main result of the
section.
\bigskip

\noindent {\bf Theorem 4.5.} {\sl Let $G$ be a finite group. Then
the following are polynomially equivalent.

(i) For every subset $A\subset G$, the directed Cayley graph
with generators in $A$ is $c_1$-quasirandom.

(ii) For every subset $A\subset G$ and every function $f:G\ra\C$
that sums to 0, $\|A*f\|\le c_2n^{1/2}|A|^{1/2}$.

(iii) Every function $f$ from $G$ to the closed unit disc in $\C$
such that $\sum_gf(g)=0$ is $c_3$-quasirandom.

(iv) For every function $f$ from $G$ to the closed unit disc in $\C$
such that $\sum_gf(g)=0$, the function $F(x,y)=f(xy^{-1})$ is
$c_3$-quasirandom on $G\times G$.

(v) Every non-trivial representation of $G$ has dimension
at least $c_4^{-1}$.}
\bigskip

\Proof  The proof that (v) implies (i) and (ii) is essentially 
contained in the argument of the previous section. Indeed, suppose
that the smallest dimension of a non-trivial representation is $k$,
and let $A\subset G$. Let $\Gamma$ be the directed Cayley graph of $A$
and let $X$ be the space of all functions $f$ such that $\sum f(x)=0$
and $\|A*f\|/\|f\|$ is maximized (together with the zero
function). Let $\lambda$ be the maximum value of this ratio. Then $X$
is invariant under the right-regular representation of $G$, so by
hypothesis it has dimension at least $k$. Lemma 2.7 implies that
$k\lambda^2\le|A|n$, so $\lambda\le(n|A|/k)^{1/2}$. This means
that if (v) holds then (ii) holds with $c_2=c_4^{1/2}$.

From this and Lemma 2.7 it follows that the number of appropriately
directed 4-cycles in $G$ is at most $|A|^4+n^2|A|^2/k$. In particular,
whatever the cardinality of $A$, the graph is at least
$k^{-1}$-quasirandom.

We proved that (iii) and (iv) were equivalent in Theorem 2.5.

Now let us prove that (iii) implies (v). That is, given a non-trivial
representation of dimension $m$, let us construct from it a function
$f$ that fails to be $c$-quasirandom for some $c$ that depends
polynomially on $m$. This we do by an averaging argument, which will
exploit the lemmas we have just proved. To simplify the notation, we
shall write the average of a function $f$ defined on the sphere $S$ as
$\E_vf(v)$ instead of $\int_Sf(v)d\mu(v)$.

A standard and easy lemma of representation theory tells us that if
$G$ has a representation $\rho$ then there is an inner product on the
vector space $V$ on which $G$ acts such that the representation is
unitary. Therefore, we may assume that $\rho$ already has this
property. Also, it will be convenient to assume, as we obviously can,
that $\rho$ is irreducible. To simplify the notation yet further, if
$v\in V$ and $g\in G$ we shall write $gv$ instead of $\rho(g)(v)$.

Given any two vectors $v$ and $w$ in the unit sphere $S$ of $V$,
let $f_{v,w}:G\ra\C$ be defined by $f_{v,w}(g)=\sp{gv,w}$. Notice
that $|f_{v,w}(g)|\le 1$ for every $g$. Furthermore, for any 
$g'$ we have
$$\sum_g gv=\sum_gg'gv=g'\Bigl(\sum_ggv\Bigr)\ .$$
Since $\rho$ is irreducible, it follows that $\sum_ggv=0$ (or
it would generate a 1-dimensional invariant subspace of $V$ and
$\rho$ would not be irreducible). Therefore,
$\sum_gf_{v,w}(g)=\sum_g\sp{gv,w}=0$. Our averaging argument
will show that at least one of these functions $f_{v,w}$ fails
to have the property in (iii), if $c_4<m^{-3}$. 

By Lemma 4.3 (for the second equality), 
$$\E_w\E_gf_{v,w}(g)\ol{f_{v,w}(gh)}=
\E_g\E_w\sp{gv,w}\sp{w,ghv}=m^{-1}\E_g\sp{gv,ghv}=m^{-1}\sp{v,hv}\ .$$
Therefore, by Lemma 4.2, 
$$\E_v\E_w\E_gf_{v,w}(g)\ol{f_{v,w}(gh)}=m^{-2}\,\ol{\tr h}\ .$$
Therefore, by the Cauchy-Schwarz inequality,
$$\E_v\E_w\Bigl|\E_gf_{v,w}(g)\ol{f_{v,w}(gh)}\Bigr|^2\ge m^{-4}|\tr h|^2\ .$$
From this it follows that
$$\E_v\E_w\E_h\Bigl|\E_gf_{v,w}(g)\ol{f_{v,w}(gh)}\Bigr|^2\ge 
m^{-4}\E_h|\tr h|^2\ ,$$
and hence that there exist $v$ and $w$ such that
$$\E_h\Bigl|\E_gf_{v,w}(g)\ol{f_{v,w}(gh)}\Bigr|^2\ge 
m^{-4}\E_h|\tr h|^2\ .$$

We now have the task of bounding $\E_h|\tr h|^2$ from below. But 
$\E_h|\tr h|^2=\E_g\E_h|\tr(gh^{-1})|^2=\E_g\E_h|\sp{A_g,A_h}|^2$,
where $A_g$ and $A_h$ are the unitary matrices corresponding to
$g$ and $h$ and the inner product comes from considering $A_g$
and $A_h$ as vectors in $\C^{m^2}$ and taking the standard 
inner product there. Since these vectors have norm $\sqrt m$, 
Lemma 4.4 implies that $\E_g\E_h|\sp{A_g,A_h}|^2\ge m$. Putting
all this together, we find that
$$\E_h\Bigl|\E_gf_{v,w}(g)\ol{f_{v,w}(gh)}\Bigr|^2\ge m^{-3}\ ,$$
completing the proof that (iii) implies (v).

All that remains to prove the theorem is to show that (i) implies
(iii). That is, given a non-quasirandom function defined on $G$, we
would like to construct from it a 01-valued function that gives
rise to a Cayley graph that is also not quasirandom. Since this
argument is standard, we shall be slightly sketchy about it.

It can be shown that the formula
$$\|F\|=\Bigl(\sum_{x,x'}\Bigl|\sum_yF(x,y)\ol{F(x',y)}\Bigr|^2\Bigr)^{1/4}$$
defines a norm $\|.\|$ on the space of functions $F:G\times G\ra\C$.
(This is a fairly easy lemma: a proof can be found in [9].) It
follows from the triangle inequality that if $F$ fails to be 
$c$-quasirandom, then either ${\rm Re} f$ or ${\rm Im} f$ fails to be 
$(c/16)$-quasirandom. Therefore, if $f$ is a function for which
(ii) fails, then there must exist a function $u$ with values in
$[-1,1]$ and average 0 such that
$$\sum_g\Bigl(\sum_hu(h)u(gh)\Bigr)^2\ge c_3|G|^3/16\ .$$
Now let $v(g)=(1+u(g))/2$ for every $g\in G$. Then a standard
argument shows that
$$\sum_g\Bigl(\sum_hv(h)v(gh)\Bigr)^2
\ge |G|^3/16+c_3|G|^3/256=(1+c_3/16)|G|^3/16\ .$$
(The argument is to expand the left-hand side into 
a sum of sixteen terms and observe that 
$$\sum_{g,g'}\Bigl(\sum_hv(h)v(gh)\Bigr)^2-
{|G|^3\over 16}-{1\over 16}\sum_{g,g'}\Bigl(\sum_hu(h)u(gh)\Bigr)^2$$
is a sum of squares.)

Now choose a subset $A\subset G$ randomly, putting $g$ into $A$ 
with probability $v(g)$, making all choices independently. Writing
$A$ also for the characteristic function of the set $A$, we wish to
estimate the sum
$$\sum_g\Bigl(\sum_hA(h)A(gh)\Bigr)^2=\sum_g\sum_{h,h'}A(h)A(gh)A(h')A(gh')
\ .$$
The number of choices of $(g,h,h')$ for which the elements $h$, $gh$,
$h'$ and $gh'$ are not all distinct is $O(|G|^2)$, and for all other
choices the expected value of $A(h)A(gh)A(h')A(gh')$ is
$v(h)v(gh)v(h')v(gh')$. Therefore, the expected value of the sum is
at least $(1+c_3/20)|G|^3/16$ when $|G|$ is sufficiently large. Also,
with very high probability $A$ has cardinality at most $(1+c_3/1000)|G|/2$
(again, if $|G|$ is sufficiently large). It follows that there exists
a set $A$ such that the directed Cayley graph defined by $A$ is not 
$c_3/32$-quasirandom.~\hfill~$\square$ \bigskip

In the light of this theorem we make the following formal definition
of a quasirandom group. Recall that quasirandom functions were defined
just after the proof of Theorem 2.5. 
\medskip

\noindent {\bf Definition.} Let $G$ be a finite group and let $c>0$.
Then $G$ is $c$-{\it quasirandom} if every function $f:G\ra\C$
that has average zero and takes values of modulus at most 1 is 
$c$-quasirandom. 

We end this section with two further characterizations of quasirandom
groups. The first one states that the quasirandom groups are precisely
those that do not contain a large product-free set. In one direction
this is the main assertion of Theorem 3.3, so we shall concentrate on
the other direction. As commented in the introduction, this final
equivalence is not a polynomial one: we shall show that if the largest
product-free subset of $G$ has size $\d|G|$, then $G$ has no
non-trivial representation of dimension less than $C\log(1/\d)$ for
some absolute constant $C$. In the final section we shall discuss
whether this result can be improved.

\proclaim Theorem {4.6}. Let $G$ be a group of order $n$ and suppose that 
$G$ has a non-trivial representation of dimension $k$. Then $G$ has a
product-free subset of size at least $c^kn$, where $c>0$ is an absolute
constant.

\Proof  Let $\phi:G\ra\C^k$ be a unitary representation of $G$. Without
loss of generality $\phi$ is irreducible, since otherwise we can find a
representation with a smaller $k$. Also, without loss of generality
it is faithful, since otherwise we can replace $G$ by
$G/\ker\phi$. Therefore, without loss of generality the elements of
$G$ are themselves unitary transformations of $\C^k$.

Now for any vector $v\in\C^k$ we have $\sum_{\a\in G}\a v=0$, since 
it is invariant under left multiplication by any $\b\in G$ and the
representation is irreducible. It follows from Lemma 4.2 that the 
average trace of an element of $G$ is 0. Since the trace of a unitary
operator has modulus at most $k$, it follows that the number of elements
$\a\in G$ such that $\tr\a$ has real part greater than $k/2$ is at most
$2n/3$. That is, at least $n/3$ elements of $G$ have trace with real
part less than or equal to $k/2$. 

Now the trace is the sum of the eigenvalues, so if $\tr\a$ has real
part at most $k/2$, there must be an eigenvalue $\omega$ with real 
part at most $1/2$. 

Let $X$ be the set of all $\a\in G$ such that $\tr\a\le k/2$ and
for each $\a\in X$ let $v(\a)$ be a unit eigenvector with eigenvalue
$\omega(\a)$ that has real part less than 1/2. 

Now let $\d>0$ be an absolute constant to be chosen later. By a
standard volume argument the unit sphere of $\C^k$ has a $\d$-net of
cardinality at most $(3/\d)^{2k}$, so we can choose at least
$(\d/3)^{2k}|X|$ elements $\a$ of $X$ such that all the vectors
$v(\a)$ lie within $\d$ of some point and hence within $2\d$ of each
other. Therefore, we can choose at least $(\d/4)^{2k}n$ elements $\a$
of $X$ such that all the $v(\a)$ are within $2\d$ of each other and
all the $\omega(\a)$ are within $\d$ of each other as well. Let $Y$ be
a subset of $X$ with this property.

We would now like to show that, for any $\a$ and $\a'$ in $Y$, the
vectors $\a v(\a)$ and $\a' v(\a)$ are close. This we deduce from the
following equalities and inequalities, which all follow from the
properties of $Y$ and the fact that the elements of $G$ preserve
distance: $\a v(\a)=\omega(\a)v(\a)$;
$\|\omega(\a)v(\a)-\omega(\a')v(\a)\|\le\d$;
$\|\omega(\a')v(\a)-\omega(\a')v(\a')\|\le 2\d$; 
$\omega(\a')v(\a')=\a'v(\a')$; $\|\a'v(\a')-\a'v(\a)\|\le 2\d$.
Therefore, by the triangle inequality, $\|\a v(\a)-\a'v(\a)\|\le 5\d$.

Now let $\a''$ be another element of $Y$. Then $\|\a
v(\a)-\a''v(\a)\|\le 5\d$ as well. Also, from the previous
inequality and the fact that $\a$ is unitary, we can deduce that
$\|\a^2v(\a)-\a\a'v(\a)\|\le 5\d$. Therefore, if $\a\a'=\a''$
it follows that $\|\a^2v(\a)-\a v(\a)\|\le 10\d$, and hence 
that $\|\a v(\a)-v(\a)\|\le 10\d$, and finally that
$|\omega(\a)-1|\le 10\d$. But we know that $\omega(\a)$ is
a complex number with modulus 1 and real part at most 1/2,
from which it follows that $|\omega(\a)-1|\ge 1$. Therefore,
$Y$ is product free as long as we choose $\d$ to be less than
$1/10$. Therefore, we can find a product-free subset $Y$ of
$G$ of size at least $c^kn$ with $c$ a positive absolute
constant (in fact, $1/2000$ will do), which proves the theorem. 
\hfill $\square$ \bigskip

Our final characterization of quasirandom groups states that a group
$G$ is quasirandom if and only if every quotient of $G$ is large and
non-Abelian. We start with a natural special case of this, showing
that {\it all non-cyclic finite simple groups are quasirandom}. One
could presumably prove this result with a better bound than we obtain
by using the classification of finite simple groups and simply looking
up the dimensions of their irreducible representations. However, our
proof is elementary. (Even this elementary argument may well be known,
but we have had trouble finding it in the literature. L\'aszlo Pyber
has pointed out to me that a slightly stronger bound can be deduced
from a theorem of Jordan, as later modified by Frobenius and
Blichfeldt, which has an elementary proof. See [10 Theorem 14.12]. 
However, the argument below is simpler.)

\proclaim Theorem {4.7}. Let $G$ be a non-cyclic finite simple 
group of order $n$. Then every non-trivial representation of $G$ 
has dimension at least $\sqrt{\log n}/2$.

\Proof  Let $\phi:G\ra U(k)$ be an irreducible unitary representation 
of $G$. Since $G$ is simple, $\phi$ has trivial kernel, so without
loss of generality $G$ itself is a finite subgroup of $U(k)$. 

Let $\a$ be any element of $G$ other than the identity. 
We claim first that $\a$ has a conjugate that does not commute with 
$\a$. To see this, suppose that all conjugates do commute with $\a$.
Then for any $\b$ and $\g$ in $G$ we have 
$$(\b\a\b^{-1})(\g\a\g^{-1})=\g(\g^{-1}\b\a\b^{-1}\g)\a\g^{-1}
=\g\a(\g^{-1}\b\a\b^{-1}\g)\g^{-1}=(\g\a\g^{-1})(\b\a\b^{-1})\ .$$
That is, all conjugates of $\a$ commute with each other. But the
subgroup of $G$ generated by conjugates of $\a$ is easily seen
to be normal, and therefore all of $G$, which implies that $G$ is 
Abelian. But in that case the only irreducible representations of 
$G$ are 1-dimensional, which implies that $k=1$ and $G$ is cyclic,
contradicting our hypothesis.

Suppose now that $\a$ is the closest element of $G$, in the
operator norm on $B(\C^k)$, to the identity (apart of course
from the identity itself), and let $\|\a-\iota\|=\e$. Let $\b$ 
be a conjugate of $\a$ that does not commute with $\a$. Then
$\|\b-\iota\|=\e$ as well, since $G$ consists of unitary
transformations. Write $\a=\iota+\g$ and $\b=\iota+\eta$. Then
$\a\b-\b\a=\g\eta-\eta\g$. Therefore, since $\a^{-1}\b^{-1}$
is unitary, $\|\iota-\a\b\a^{-1}\b^{-1}\|=\|\g\eta-\eta\g\|$.
Since $\a$ and $\b$ do not commute, and are closest
elements to the identity, it follows that $\|\g\eta-\eta\g\|\ge\e$.
But we also know that $\|\g\eta-\eta\g\|\le 2\|\g\|\|\eta\|=2\e^2$.
Therefore, $\e\ge 1/2$, which implies that no two elements of
$G$ are closer than $1/2$ in the operator norm.

It remains to determine an upper bound for the size of a 
$1/2$-separated subset of $U(k)$. But $U(k)$ is contained in
the unit ball of $B(\C^k)$. The volume argument mentioned
in the previous lemma shows that for any $d$-dimensional 
real normed space and any $\e>0$ the largest $\e$-separated
subset of the unit ball has size at most $(1+2/\e)^d$. 
The normed space $B(\C^k)$ is a $k^2$-dimensional complex
space, so, setting $d=2k^2$ and $\e=1/2$, we deduce that a
$1/2$-separated subset of $U(k)$ has cardinality at most
$25^{k^2}$. That is, $n\le 25^{k^2}$, from which the 
theorem follows. \hfill $\square$ \bigskip

Note that the alternating groups $A_n$ have representations of
dimension $n-1$ (since they act on the subspace of $\C^n$ consisting
of vectors whose coordinates add up to 0). Therefore, the bound in
Theorem 4.7 cannot be improved to more than $\log n/\log\log n$.



\proclaim Theorem {4.8}. Let $G$ be a group of order $n$
and suppose that for every proper normal subgroup $H$ of $G$, the
quotient $G/H$ is non-Abelian and has order at least $m$. Then $G$ has
no non-trivial representation of dimension less than $\sqrt{\log m}/2$.
Conversely, if $G$ has an Abelian quotient, then $G$ has a 1-dimensional
representation, and if $G$ has a quotient of order $m$, then $G$ has
a representation of dimension $\sqrt m$.

\Proof  Let us quickly deal with the converse, since this is easy 
and not the main point of interest. Any representation of a quotient
of $G$ can be composed with the quotient map so that it becomes a
representation of $G$ of the same dimension. Therefore, the result
follows from two standard facts of representation theory: that 
the irreducible representations of Abelian groups are 1-dimensional
(and exist!), and that every group of order $m$ has a representation
of dimension at most $\sqrt m$. (This second fact follows from the
result that the sum of the squares of the dimensions of the irreducible
representations is $m$.)

Now let us turn to the more interesting direction of the theorem.
Let $H$ be a maximal proper normal subgroup of $G$. Then the 
quotient group $G/H$ is simple and, by our hypothesis, non-Abelian.
Let $\phi:G\ra U(k)$ be a unitary representation of $G$. 
If we knew that the kernel of $\phi$ was $H$, then we would have
a representation of $G/H$ to which we could apply Theorem 4.7. 
However, this does not have to be the case, so instead we must
imitate the proof of Theorem 4.7, as follows.

We may clearly assume that $\phi$ is a faithful representation
(or else we look at the quotient of $G$ by its kernel). Therefore,
we shall think of the elements of $G$ itself as unitary maps
on $\C^k$. Let us now define a metric on $G/H$ by taking 
$d(\a H,\b H)$ to be the smallest distance (in the operator 
norm again) between any element of $\a H$ and any
element of $\b H$. Let $\a$ be an element of $G\setminus H$ 
such that the distance from $\a H$ to $H$, with respect to this
metric, is minimized, and note that this distance is just the
smallest distance in the operator norm from any element of 
$\a H$ to the identity. Without loss of generality, $\a$ itself
is an element of $\a H$ for which this minimum is attained.

Now $G/H$ is simple and non-Abelian. Hence, by the argument of
the last section, we can find a conjugate $\b H$ of $\a H$ 
in $G/H$ that does not commute with $\a H$. It is easy to 
see that we can choose the representative $\b$ to be a 
conjugate of $\a$ in $G$, so let us do this. Then $\b$ is a conjugate
of $\a$ such that not only do $\a$ and $\b$ not commute, but
they do not even belong to the same coset of $H$. Moreover,
the distance from $\b$ to the identity is the same as the
distance from $\a$ to the identity. As in the proof of Theorem 
4.8, let $\e$ be this distance, and let $\a=\iota+\gamma$ and
$\b=\iota+\eta$. 

Once again, the distance between $\a\b$ and $\b\a$ is 
$\|\g\eta-\eta\g\|$, and therefore so is the distance between
$\iota$ and $\a\b\a^{-1}\b^{-1}$. Since $\a\b\a^{-1}\b^{-1}$
does not belong to $H$, it follows from our minimality assumption
that $\|\g\eta-\eta\g\|\ge\e$, as before, and it is also at 
most $2\e^2$ for precisely the same reason as before. Therefore,
no two elements of different cosets of $H$ can be within $1/2$
of each other in the operator norm, so, by the upper bound given
in the proof of Theorem 4.7 for the size of a $1/2$-separated 
subset of $U(k)$, there can be at most $25^{k^2}$ cosets of $H$.
This proves the theorem. \hfill $\square$ \bigskip

A good example to bear in mind in connection with Theorem 4.8 and its
proof is the following family of groups. Let $p$ and $k$ be positive
integers and let $G(p,k)$ be the subgroup of $U(k)$ generated by all
diagonal matrices with $p$th roots of unity as their diagonal entries,
and all permutation matrices corresponding to even permutations. Thus,
a typical element of $G(p,k)$ is a permutation matrix of determinant 1
with its 1s replaced by arbitrary $p$th roots of unity. The subgroup
$H(p,k)$ generated by just the diagonal matrices in $G(p,k)$ is
normal, and the quotient is isomorphic to the alternating group
$A_k$. Moreover, one can show that any proper normal subgroup of
$G(p,k)$ is contained in $H(p,k)$. Therefore, these groups are
quasirandom as $k$ tends to infinity, despite being of arbitrarily
high order for any fixed $k$. The reason this can happen is that, as
the proof of Theorem 4.8 shows is necessary, the cosets of $H(p,k)$
are well-separated.


\bigskip


In practice, Theorems 4.6 and 4.8 are not particularly useful
characterizations of quasirandomness because the equivalences are not
polynomial equivalences. In other words, they are fine if all one
wants is qualitative statements (such as that no subset of positive
density is product free) but too crude if one is interested in bounds
of the kind obtained in this paper. However, sometimes a qualitative
statement is interesting -- for example, if one is wondering whether
a particular family of groups is quasirandom and wants to make a 
preliminary check. For instance, Theorem 4.8 tells us that SL$_2(p)$ 
is quasirandom, since $\{\iota,-\iota\}$ is a maximal normal subgroup 
of very high index. However, this particular group is much more 
quasirandom than Theorem 4.8 guarantees. As for Theorem 4.6, it can
in fact be improved to a polynomial equivalence: this will be discussed
in the final section.
\bigskip

\noindent {\bf \S 5. Solving equations in quasirandom groups.}
\bigskip

The purpose of this section is to prove a generalization of
Theorem 3.3: instead of finding $a$ and $b$ such that $a$, $b$ 
and $ab$ each lie in specified sets, we shall find $\seq a m$
such that for every non-empty subset $F\subset\{1,2,\dots,m\}$ 
the product of those $a_i$ with $i\in F$ lies in a specified set.
In other words, perhaps surprisingly, we can choose $m$ elements
of the group in such a way that exponentially many conditions are 
satisfied simultaneously, using only the fact that a reasonable
number of elements satisfy each condition individually.

Underlying the argument is the following basic lemma, which is
a reformulation of the last part of Theorem 3.3 that will be
slightly more convenient. The proof of the main theorem
of this section will use it to drive an inductive argument.

\proclaim Lemma {5.1}. Let $G$ be a group of order $n$ such
that no non-trivial representation has dimension less than $k$.
Let $A$ and $B$ be two subsets of $G$ with densities $rn$ and
$sn$, respectively and let $\d$ and $t$ be two positive
constants. Then, provided that $rst\ge(\d^2k)^{-1}$, 
the number of group elements $x\in G$ for which 
$|A\cap xB|\le(1-\d)rsn$ is at most $t n$.

\Proof  Let $C$ be the set $\{x^{-1}:x\in B\}$. Then
$$|A\cap xB|=\sum_yA(y)(xB)(y)=\sum_yA(y)B(x^{-1}y)=\sum_yA(y)C(y^{-1}x)
=A*C(x)\ .$$
By Theorem 4.5, if $f:G\ra\R$
sums to zero, then $\|A*f\|\le(r/k)^{1/2}n\|f\|$.
Applying this result in the case $f(x)=C(x)-s$ and noting
that $\|f\|^2=s(1-s)n\le sn$, we deduce
that $\|A*C-rsn\|^2\le rsn^3/k$. It follows that the number
of $x$ such that $A*C(x)\le(1-\d)rsn$ is at most 
$n/\d^2rsk$. If $rst\ge(\d^2k)^{-1}$, then
this is at most $tn$, as required. \hfill $\square$ \bigskip 

Note the following easy consequence of Lemma 5.1, which shows that it
is indeed effectively the same as Theorem 3.3. Suppose that $rst>1/k$
and that $C$ is a subset of $G$ with density $t$.  Lemma 5.1 with
$\d=1$ tells us that the number of $y$ such that $A\cap
y^{-1}B=\emptyset$ is less than $tn$, from which it follows that there
exists $y\in C$ such that $A\cap y^{-1}B\ne\emptyset$.  But then, if
$x\in A\cap y^{-1}B$, we have $x\in A$, $y\in C$ and $yx\in B$.

In order to make the proof of our general theorem more transparent, 
we begin with the special case $m=3$.

\proclaim Theorem {5.2}. Let $G$ be a group of order $n$ such
that no non-trivial representation has dimension less than $k$.  Let
$A_1$, $A_2$, $A_3$, $A_{12}$, $A_{13}$, $A_{23}$ and $A_{123}$ be
subsets of $G$ of densities $p_1$, $p_2$, $p_3$, $p_{12}$, $p_{13}$,
$p_{23}$ and $p_{123}$, respectively. Then, provided that $p_1p_2p_{12}$,
$p_1p_3p_{13}$, $p_1p_{23}p_{123}$ and $p_2p_3p_{23}p_{12}p_{13}p_{123}$ are
all at least $16/k$, there exist elements $x_1\in A_1$, $x_2\in A_2$
and $x_3\in A_3$ such that $x_1x_2\in A_{12}$, $x_1x_3\in A_{13}$,
$x_2x_3\in A_{23}$ and $x_1x_2x_3\in A_{123}$.

\Proof  We start by choosing $x_1$, noting that there are certain
conditions it will have to satisfy if there is to be any hope of
continuing the proof. For example, later we shall need to choose
$x_2\in A_2$ such that $x_1x_2\in A_{12}$. Equivalently, we shall
need $x_2$ to belong to $A_2\cap x_1^{-1}A_{12}$. Similarly, we
shall need $x_3\in A_3\cap x_2^{-1}A_{13}$ and 
$x_2x_3\in A_{23}\cap x_1^{-1}A_{123}$. Therefore, we want these
sets to be not just non-empty, but reasonably large.

By Lemma 5.1, the number of $x_1$ such that 
$|A_2\cap x_1^{-1}A_{12}|< p_2p_{12}n/2$ is at most $p_1n/4$, 
provided that $p_1p_2p_{12}\ge 16/k$. Similarly, if 
$p_1p_3p_{13}\ge 16/k$ and $p_1p_{23}p_{123}\ge 16/k$, then
the number of $x_1$ such that $|A_3\cap x_1^{-1}A_{13}|<p_3p_{13}n/2$
is at most $p_1n/4$ and the number of $x_1$ such that
$|A_{23}\cap x_1^{-1}A_{123}|< p_{23}p_{123}n/2$ is at most $p_1n/4$.
Therefore, provided these inequalities hold, we can choose $x_1\in A_1$
such that, setting $B_2=A_2\cap A_{12}$, $B_3=A_3\cap A_{13}$ and
$B_{23}=A_{23}\cap A_{123}$, $q_2=p_2p_{12}/2$, $q_3=p_3p_{13}/2$
and $q_{23}=p_{23}p_{123}/2$, we have $|B_2|\ge q_2n$, $|B_3|\ge q_3n$
and $|B_{23}|\ge q_{23}n$. 

At this point we could quote our results about product-free sets, but
instead let us repeat the argument (which is more or less an equivalent
thing to do). We would like to choose $x_2\in B_2$ such that 
$B_3\cap x_2^{-1}B_{23}$ is non-empty. Lemma 5.1 implies that
the number of $x_2$ such that $B_3\cap x_2^{-1}B_3$ is empty
is at most $q_2n/2$, provided that $q_2q_3q_{23}\ge 2/k$.
Therefore, provided we have this inequality, which, when expanded,
says that $p_2p_3p_{23}p_{12}p_{13}p_{123}\ge 16/k$, there exist
$x_2\in B_2$ and $x_3\in B_3$ such that $x_2x_3\in B_{23}$. But
then $x_1$, $x_2$ and $x_3$ satisfy the conclusion of the theorem.
\hfill $\square$ \bigskip

It is clear that the above argument can be generalized. The only thing
that is not quite obvious is the density conditions that emerge from the
resulting inductive argument. Here is what they are. Suppose
that for every subset $F\subset\{1,2,\dots,m\}$ we have a subset $A_F$
of a group $G$ with density $p_F$ and suppose that no non-trivial
representation of $G$ has dimension less than $k$. Now let $h$ be an
integer less than $m$ and let $E$ be a subset of $\{h+1,\dots,m\}$. 
Let $\ca_{h,E}$ be the collection of all sets of the form $U\cup V$, 
where $\max U<h$ and $V$ is either $\{h\}$, $E$ or $\{h\}\cup E$.
We shall say that the sets $A_F$ {\it satisfy the $(h,E)$-density 
condition} if $\prod_{F\in\ca_{h,E}}p_F$ is at least $2^{3m}/k$. We 
shall say that they {\it satisfy the density condition} if they 
satisfy the $E$-density condition for every $h<m$ and every
non-empty set $E\subset\{h+1,\dots,m\}$.

To get an idea of what this means, notice that the inequalities we
assumed in Theorem 5.2 are the $(1,\{2\})$-condition, the
$(1,\{3\})$-condition, the $(1,\{2,3\})$-condition and the
$(2,\{3\})$-condition, respectively, except that there we had
a slightly better dependence on $m$.

\proclaim Theorem {5.3}. Let $G$ be a group of order $n$ such
that no non-trivial representation has dimension less than $k$. 
For each non-empty subset $F\subset\{1,2,\dots,m\}$ let $A_F$
be a subset of $G$ of density $p_F$, and suppose that this 
collection of sets satisfies the density condition. Then there
exist elements $\seq x m$ of $G$ such that $x_F\in A_F$ for
every $F$, where $x_F$ stands for the product of all $x_i$
such that $i\in F$, written with the indices in increasing order.

\Proof  By the density condition, for every non-empty subset 
$F\subset\{2,\dots,m\}$ we have the inequality $2^{-m}p_1p_Fp_{1F}\ge
2^{2m}/k$. (Here we use the shorthand $1F$ to stand for $\{1\}\cup F$.)
Therefore, by Lemma 5.1, for each $F$ the number of $x_1$ such that
$|A_F\cap x^{-1}A_{1F}|\le p_Fp_{1F}(1-2^{-m})$ is at most
$p_1n/2^m$. Therefore, the number of $x_1$ such that $|A_F\cap
x_1^{-1}A_{1F}|\le p_Fp_{1F}(1-2^{-m})$ for at least one non-empty
$F\subset\{2,\dots,m\}$ is at most $p_1n/2$. It follows
that there exists $x_1\in A_1$ such that, if for every non-empty
$F\subset\{2,\dots,m\}$ we set $B_F=A_F\cap A_{1F}$, then every $B_F$
has density at least $q_F=p_Fp_{1F}(1-2^{-m})$.

We claim now that the sets $B_F$ satisfy the density condition (after
a relabelling of the index set). Let $h<m$ and let $E$ be a non-empty
subset of $\{h+1,\dots,m\}$. Define $\cb_{h,E}$ to be the set of all
$F$ of the form $U\cup V$ with $U\subset\{2,\dots,h-1\}$ and $V$ equal
to $\{h\}$, $E$ or $\{h\}\cup E$. Then
$$\prod_{F\in\cb_{h,E}}q_F
\ge(1-2^{-m})^{2^m}\prod_{F\in\cb_{h,E}}p_Fp_{1F}
=(1-2^{-m})^{2^m}\prod_{F\in\ca_{h,E}}p_F\ .$$
But $(1-2^{-m})^{2^m}\ge 1/4$ and $\prod_{F\in\ca_{h,E}}p_F\ge 2^{3m}/k$,
so this implies that $\prod_{F\in\cb_{h,E}}q_F\ge 2^{3(m-1)}/k$.
Therefore, the sets $B_F$ satisfy the density condition.

This proves the inductive step of the theorem. To be on the safe side,
we take as our base case the case $m=2$. (We do this so that we do 
not have to worry about the definition of the density condition when
$E$ cannot be non-empty.) This follows easily from the remark following
Lemma 5.1 if one sets $A_1=C$, $A_2=B$ and $A_{12}=A$. The density 
condition in this case is stronger than the hypothesis we needed 
to guarantee the existence of $x_1$ and $x_2$ such that $x_1\in A_1$,
$x_2\in A_2$ and $x_{12}\in A_{12}$. Therefore, the theorem is 
proved.~\hfill~$\square$~\bigskip

We now give a couple of corollaries of Theorem 5.3. They are special
cases of the theorem: the only extra content is that we need to do a
small amount of calculation to optimize certain densities while 
preserving the density condition.

\proclaim Corollary {5.4}. Let $G$ be a group of order $n$ such
that no non-trivial representation has dimension less than $k$. 
For each non-empty subset $F\subset\{1,2,\dots,m\}$ let $A_F$
be a subset of $G$ of density $p$. Then, provided that 
$p^{3.2^{m-2}}>2^{3m}/k$ (which is true if $p>2k^{-1/2^{2m}}$),
there exist $\seq x m$ such that $x_F\in A_F$ for every $F$.

\Proof Since all the densities are the same, all we have to do
is look at which set $\ca_{h,E}$ is largest. Obviously they get
larger as $h$ gets larger, so the largest one is when $h=m-1$.
This has size $3.2^{m-2}$ since there are $2^{m-2}$ possibilities
for $U$ and 3 possibilities for $V$. The result now follows from
Theorem 5.3. \hfill $\square$ \bigskip

\proclaim Corollary {5.5}. Let $G$ be a group of order $n$ such
that no non-trivial representation has dimension less than $k$. 
For every pair $1\le i<j\le m$ let $A_{ij}$ be a set of density
$p$. Then, provided that $p>4k^{-1/(2m-3)}$, there exist $\seq x m$ 
such that $x_ix_j\in A_{ij}$ for every $i<j$.

\Proof  We shall apply Theorem 5.3 again, setting $A_F$ to be
$G$ whenever $F$ has cardinality other than 2. Then $p_F=p$
if $F$ has cardinality 2, and $p_F=1$ otherwise. Now let us
work out how many sets of size 2 are contained in $\ca_{h,E}$.
If $E$ has cardinality greater than 1 then there are $h-1$
such sets, since then $V$ must equal $\{h\}$ and $U$ must be a 
singleton. If $E$ has cardinality equal to 1 then there are
$2h-1$ sets, since either $U$ is a singleton and $V$ is $\{h\}$
or $E$, or $U$ is empty and $V$ is $\{h\}\cup E$. Since the
largest possible value of $h$ is $m-1$, this tells us that
the sequence exists provided that $p^{2m-3}>2^{3m}/k$, which
implies the corollary. \hfill $\square$ \bigskip

It is possible to generalize Theorem 5.3 slightly further by
exploiting two facts about Lemma 5.1. Instead of giving full details,
we shall merely state two results and briefly explain how they are
proved.

\proclaim Theorem {5.6}. Let $G$ be a group of order $n$ such
that no non-trivial representation has dimension less than $k$.  
For every pair $1\le i<j\le m$ let $A_{ij}$ be a set of density
$p$. Then, provided that $p>4k^{-1/(2m-3)}$, there exist $\seq x m$ 
such that $x_ix_j^{-1}\in A_{ij}$ for every $i<j$.

\proclaim Theorem {5.7}. Let $G$ be a group of order $n$ such
that no non-trivial representation has dimension less than $k$.
Let $A_1$, $A_2$, $A_3$, $A_{12}$, $A_{13}$, $A_{23}$ and $A_{123}$ be
subsets of $G$ of densities $p_1$, $p_2$, $p_3$, $p_{12}$, $p_{13}$,
$p_{23}$ and $p_{123}$, respectively. Then, provided that $p_1p_2p_{12}$,
$p_1p_3p_{13}$, $p_1p_{23}p_{123}$ and $p_2p_3p_{23}p_{12}p_{13}p_{123}$ are
all at least $16/k$, there exist elements $x_1\in A_1$, $x_2\in A_2$
and $x_3\in A_3$ such that $x_1x_2\in A_{12}$, $x_3x_1\in A_{13}$,
$x_2x_3^{-1}\in A_{23}$ and $x_2x_3^{-1}x_1^{-1}\in A_{123}$.

To prove statements like this, one exploits Lemma 5.1 and its
method of proof to the full. Not only can one show that $A\cap xB$
is nearly always about the same size (when $A$ and $B$ are large
enough), but also $A\cap x^{-1}B$, $A\cap Bx$ and $A\cap Bx^{-1}$. 
The inductive proof of Theorem 5.3 works as long as at each stage
of the inductive process the variable one is trying to choose,
or its inverse, appears either at the beginning or at the end of 
each product. So, for example, in Theorem 5.7 one starts by choosing
$x_1$ such that $A_2\cap x_1^{-1}A_{12}$, $A_3\cap A_{13}x_1^{-1}$
and $A_{23}\cap A_{123}x_1$ are all large. One is then left needing
to place $x_2$, $x_3$ and $x_2x_3^{-1}$ into these sets, which can
clearly be done.
\bigskip

\noindent {\bf Remarks.} Although it may at first seem surprising that 
one can cause so many equations to be satisfied simultaneously, there
is an intuitive explanation for this, at least for readers familiar
with the notion of higher-degree uniformity for subsets of Abelian
groups. (See [8, Section 3] for a definition of this.)  In that
terminology, Lemma 5.1 shows that all dense subsets of $G$ have a
property very similar to uniformity. But if that is the case, then
almost all intersections of a dense set $A$ with a translate of itself
will still be dense, and will therefore be uniform as well, which
shows that $A$ has a sort of non-Abelian version of quadratic
uniformity. But if uniformity implies quadratic uniformity, then it
implies uniformity of all degrees.  In the Abelian case, the higher
the degree of uniformity a set has, the more linear equations one can
hope to solve simultaneously in that set, so it is not too surprising
after all that one can solve large numbers of equations simultaneously
in subsets of a group where every dense set is uniform.

Another interesting aspect of Theorem 5.3 is that under certain 
circumstances it can yield very good bounds. For simplicity let
us consider the case where all the sets $A_F$ have density
either $p$ or $1$, and let $\cal F$ be the set of $F$ such that
the density is $p$. Suppose that no element of $\{1,2,\dots,m\}$
is contained in more than $r$ of the sets $F\in\cal F$. Then 
no set $\ca_{h,E}$ can contain more than $2r$ elements of $\cal F$,
so we can satisfy all the conditions simultaneously if 
$p^{2r}\ge 2^{3m}/k$. That is, for fixed $r$ we can contain a
power that is independent of $m$. (With a bit of care, the
exponential dependence of the constant on $m$ can be improved
as well.) This situation would arise if, for example, we wanted
$x_ix_j$ to belong to $A_{ij}$ whenever $ij$ was an edge of a
certain graph $H$ of maximal degree 10.
\bigskip

\noindent {\bf \S 6. Open questions.}
\bigskip

The results of this paper leave several questions unanswered. One that
has been mentioned already is the following (which is not formulated
in a precise manner).

\proclaim Question {6.1}. Is there a good model for large random finite
groups with the property that a group chosen according to this model
has a high probability of being quasirandom?

Another question that has been touched on is whether Theorem
4.6 can be improved. More precisely, in an earlier draft of this
paper the following was asked. 

\proclaim Question {6.2}. If $G$ has a non-trivial
representation of dimension $k$, does $G$ have a product-free subset
of size $cn$ for some $c$ that depends polynomially on $k^{-1}$? 

I am grateful to L\'azslo Pyber for informing me that the answer is
yes, for the following reason. It can be shown using the
classification of finite simple groups that a finite group with a
$k$-dimensional representation must have a proper subgroup of index at
most $k^c$ (for some absolute constant $c$) or an Abelian quotient.
But in both cases it is easy to construct product-free subsets. A
stronger result that also implies a positive answer to Question 6.2
can be found in a recent paper of Nikolov and Pyber [15]. This 
leaves open the question of whether the classification of finite
simple groups is needed for solving Question 6.2. The results
used in the solution just mentioned do seem to have that flavour,
but it does not seem completely unreasonable to hope for a 
classification-free answer to the question. We put this as our
next question.

\proclaim Question {6.3}. Is there an elementary proof that if
$G$ has a non-trivial representation of dimension $k$ then $G$
has a product-free subset of size $cn$ for some $c$ that depends
polynomially on $k^{-1}$?

A closely related question is to find good bounds for the
largest Haar measure of a product-free subset of $SU(n)$. The methods
of this paper, suitably adapted, ought to prove that this is at most
$Cn^{-1/3}$, but the largest product-free subsets of $SU(n)$ that
we know of are in the spirit of the construction of Theorem 4.6 and
are therefore exponentially small. We therefore ask the following
question, with a tentative expectation that the answer is yes.

\proclaim Question {6.4}. Does there exist a constant $c<1$ such
that every subset $A\subset SU(n)$ that is measurable and product-free
has measure at most $c^n$?

It is easy to prove that no stronger bound can hold: just fix a unit
vector $x_0\in\C^n$ and let $A$ be the set of unitary maps $\alpha$
such that $\sp{x_0,\alpha x_0}<-1/2$. If $\alpha$, $\beta$ and $\alpha\beta$
all belong to $A$, then $\sp{x_0,\alpha x_0}$, $\sp{x_0,\alpha\beta x_0}$
$\sp{\alpha x_0,\alpha\beta x_0}$ are all less than $-1/2$. But it is
an easy exercise to show that it is impossible to find three unit
vectors with this property. (Just look at the square of the norm of 
their sum.) It is also easy to see that $A$ has size at least $c^n$
for some positive constant $c$.

Several problems arise when one starts to think about the following
broad question: which equations have solutions in large subsets of
PSL$_2(q)$, or of other quasirandom groups?  The most general answer
we have been able to find is Theorem 5.3 (and the slight generalization
mentioned at the end of the last section), but it is not obvious that
that is the end of the story. Here are two questions that give some
idea of what further results might or might not be true. The first has
an easy negative answer: if $A$, $B$ and $C$ are three large sets, can
one find $a\in A$, $b\in B$ and $c\in C$ such that $ab=ca$? The answer
is no, since if $ab=ca$, then $b=a^{-1}ca$. Thus, $b$ and $c$ are
conjugate, so to find a counterexample all one has to do is
make $B$ and $C$ disjoint unions of conjugacy classes.

However, for a very similar question it is much less clear what
the answer is. If $A$ is a quasirandom subset of an Abelian group,
then $A$ contains approximately the same number of arithmetic 
progressions of length 3 (defined to be sequences of the form 
$(a,a+d,a+2d)$ with $d\ne 0$) as a random set of the same 
cardinality, and it also contains about the same number of
solutions to the equation $x+y=z$. Moreover, the proofs of these
two facts are very similar. What happens if we investigate
arithmetic progressions in subsets of PSL$_2(q)$? 

The most obvious question is not very interesting: does every dense
subset $A$ of PSL$_2(q)$ contain a progression of length 3, where this
is now defined to be a sequence of the form $(x,gx,g^2x)$? (It might
be better to call this a ``left progression,'' since it is not the
same as a sequence of the form $(x,xg,xg^2)$.) The answer is yes,
since PSL$_2(q)$ can be decomposed into right cosets of a cyclic
subgroup of order $q$: we can therefore find a coset such that $A$
intersects it densely and apply Roth's theorem.  However, this leaves
two questions unanswered. The first is whether $A$ must in fact 
contain roughly the ``expected'' number of progressions of length 3.

\proclaim Question {6.5}. Let $A$ be a subset of PSL$_2(q)$ of
density $\d$ and let $g$ and $x$ be randomly chosen elements of
PSL$_2(q)$. Is the probability that $x$, $gx$ and $g^2x$ are
all in $A$ necessarily approximately equal to $\d^3$? 

The second question is closely related.

\proclaim Question {6.6}. Let $A$, $B$ and $C$ be three dense subsets 
of PSL$_2(q)$. Must there be an arithmetic progression 
$(a,b,c)\in A\times B\times C$?

This would be interesting, since an ``off-diagonal'' Roth theorem of
this kind is completely false in an Abelian group. Of course, the
last two questions can be asked for other quasirandom groups.
Notice also that if $(a,b,c)=(x,gx,g^2x)$, then $c=ba^{-1}b$, and if
$c=ba^{-1}b$ then $(a,b,c)=(a,ga,g^2a)$ for $g=ba^{-1}$. Therefore,
an equivalent question to the last one is the following: if $A$, $B$ and
$C$ are three dense subsets of PSL$_2(q)$, must there exist $a\in A$,
$b\in B$ and $c\in C$ such that $bab=c$? (To make the question cleaner
we have replaced $A$ by the set of inverses of elements of $A$, which
obviously makes no difference.)

There is a natural bipartite graph that one can define in response to
these problems: join $x$ to $y$ if there exists $b\in B$ such that
$bxb=y$. If this graph is automatically quasirandom, then the answers
to both problems are yes. But it is not clear whether it is
quasirandom. The difficulty is that we are mixing left and right
actions, which makes representation theory less easy to apply.
(Notice that the natural bipartite graph associated with the equation
$ab=ca$ we considered first joins $x$ to all points of the form
$a^{-1}xa$.  It is easy to see that this graph is very far from
quasirandom -- indeed, it has multiple edges and a typical edge has
very high multiplicity.)
\bigskip

\noindent {\bf Acknowledgements.} I am grateful to Vera S\'os for
drawing my attention to this problem, and to L\'aszl\'o Babai, Alexander
Gamburd, Kiran Kedlaya, L\'aszl\'o Pyber, Vlado Nikiforov and the
referee for useful remarks, especially concerning certain parts of the paper
where I trespassed into areas about which I knew very little.
\bigskip

\noindent {\bf References.}
\bigskip

\noindent [1] L. Babai, L. R\'onyai, {\it Computing irreducible
representations of finite groups}, Math. Comp. {\bf 55} (1990),
705-722.
\medskip

\noindent [2] L. Babai, V. S\'os, {\it Sidon sets in groups and 
induced subgraphs of Cayley graphs}, European J. Combin. {\bf 6}
(1985), 101-114.
\medskip

\noindent [3] B. Bollob\'as and V. Nikiforov, {\it Hermitian matrices
and graphs: singular values and discrepancy}, Discrete Math. {\bf 285}
(2004), 17-32.
\medskip

\noindent [4] J. Bourgain, A. Gamburd, {\it Uniform expansion bounds
for Cayley graphs of SL$_2(F_p)$}, preprint.
\medskip

\noindent [5] F.R.K. Chung, R.L. Graham, {\it Quasi-random subsets of 
$\Bbb Z_n$}, J. Comb. Th. A {\bf 61} (1992), 64--86.
\medskip

\noindent [6] F.R.K. Chung, R.L. Graham, R.M. Wilson, {\it Quasi-random 
graphs}, Combinatorica {\bf 9} (1989), 345--362.
\medskip

\noindent [7] G. Davidoff, P. Sarnak, A. Valette, Elementary number
theory, group theory, and Ramanujan graphs, London Mathematical Society
Student Texts, 55, Cambridge University Press, Cambridge, 2003.
\medskip

\noindent [8] W. T. Gowers, {\it A new proof of Szemer\'edi's theorem},
Geom. Funct. Anal. {\bf 11} (2001), 465-588.
\medskip

\noindent [9] W. T. Gowers, {\it Quasirandomness, counting and regularity
for 3-uniform hypergraphs}, Combin. Probab. Comput. {\bf 15} (2006),
143-184.
\medskip

\noindent [10] I. M. Isaacs, Character theory of finite groups, AMS 
Chelsea Publishing, Providence, RI, 2006 (corrected reprint of 1976
original), xii+310pp.
\medskip

\noindent [11] K. S. Kedlaya, {\it Large product-free subsets of finite
groups}, J. Combin. Theory Ser. A {\bf 77} (1997), 339-343.
\medskip

\noindent [12] K. S. Kedlaya, {\it Product-free subsets of groups},
Amer. Math. Monthly {\bf 105} (1998), 900-906.
\medskip

\noindent [13] K. S. Kedlaya, {\it Product-free subsets of groups, then 
and now}, preprint, arXiv:0708.2295v1.
\medskip

\noindent [14] A. Lubotzky, R. Phillips, P. Sarnak, {\it Ramanujan
graphs}, Combinatorica {\bf 8} (1988), 261-277.
\medskip

\noindent [15] N. Nikolov, L. Pyber, {\it Product decompositions of quasirandom
groups and a Jordan type theorem}, preprint, arXiv:math/0703343v3.
\medskip

\noindent [16] P. Sarnak and X. Xue, {\it Bounds for multiplicities of
automorphic representations}, Duke Math. J. {\bf 64} (1991), 207-227.
\medskip

\noindent [17] A. G. Thomason, {\it Pseudo-random graphs}, Proceedings 
of Random Graphs, Pozn\'an 1985 (M. Karonski, ed.), Annals of Discrete 
Mathematics {\bf 33}, 307--331.
\medskip

\bye